\documentclass{article}

\usepackage[preprint]{neurips_2023_without_decoration}

\usepackage[utf8]{inputenc}  %
\usepackage[T1]{fontenc}     %
\usepackage{hyperref}        %
\usepackage{url}             %
\usepackage{booktabs}        %
\usepackage{amsfonts}        %
\usepackage{nicefrac}        %
\usepackage{microtype}       %
\usepackage{xcolor}          %
\title{The Minimax Rate of HSIC Estimation for Translation-Invariant Kernels}

\author{%
  Florian Kalinke \\
  Institute for Program Structures and Data Organization\\
  Karlsruhe Institute of Technology\\
  Karlsruhe, Germany\\
  \texttt{florian.kalinke@kit.edu} \\
  \And
  Zoltán Szabó \\
  Department of Statistics\\
  London School of Economics\\
  London, UK \\
  \texttt{z.szabo@lse.ac.uk}
}

\usepackage{mathtools}
\mathtoolsset{showonlyrefs}  
\usepackage{amsmath}
\usepackage{amsthm}
\usepackage{bm}
\usepackage{enumitem}

\usepackage[mathscr]{eucal}

\newcommand*{\T}{^{\mkern-1.5mu\mathsf{T}}} %
\newcommand{\E}{\mathbb{E}}   %
\newcommand{\Q}{\mathbb{Q}}   %
\newcommand{\R}{\mathbb{R}}   %
\newcommand{\X}{\mathcal{X}}  %
\newcommand{\ip}[2]{\left\langle{#1}\right\rangle_{#2}} %
\newcommand{\norm}[2]{\left\|{#1}\right\|_{#2}} %
\newcommand{\tb}{\textbf}    %
\renewcommand{\H}{\mathcal{H}} %
\renewcommand{\O}{\mathcal{O}} %
\renewcommand{\P}{\mathbb{P}} %
\renewcommand{\Q}{\mathbb{Q}} %
\renewcommand{\b}{\mathbf}    %
\renewcommand{\d}{\mathrm{d}} %
\newcommand{\bdiag}{\operatorname{bdiag}} %

\DeclareMathOperator{\trace}{tr}

\DeclareMathOperator{\HSIC}{HSIC}
\DeclareMathOperator{\MMD}{MMD}

\setenumerate{labelindent=0em,leftmargin=1.6em,topsep=0.1cm,partopsep=0cm,parsep=0cm,itemsep=2mm}
\setitemize{labelindent=0em,leftmargin=1.3em,topsep=0.1cm,partopsep=0cm,parsep=0cm,itemsep=2mm}

\newtheorem*{theorem*}{Theorem}
\newtheorem{theorem}{Theorem}
\newtheorem{lemma}{Lemma}
\newtheorem{corollary}{Corollary}
\newtheorem{remark}{Remark}

\newtheorem{theoremA}{Theorem}[section]
\newtheorem{lemmaA}{Lemma}[section]

\begin{document}

\maketitle

\begin{abstract}
  Kernel techniques are among the most influential approaches in data science and statistics. Under mild conditions, the reproducing kernel Hilbert space associated to a kernel is capable of encoding the independence of $M\ge2$ random variables.
  Probably the most widespread independence measure relying on kernels is the so-called Hilbert-Schmidt independence criterion (HSIC; also referred to as distance covariance in the statistics literature). Despite various existing HSIC estimators designed since its introduction close to two decades ago, the fundamental question of the rate at which HSIC can be estimated is still open.
  In this work, we prove that the minimax optimal rate of HSIC estimation on $\R^d$ for Borel measures containing the Gaussians with continuous bounded translation-invariant  characteristic  kernels  is $\O\!\left(n^{-1/2}\right)$. Specifically, our result implies the optimality in the minimax sense of many of the most-frequently used estimators (including the U-statistic, the V-statistic, and the Nyström-based one) on $\R^d$.
\end{abstract}

\section{Introduction}

Kernel methods
\citep{steinwart08support,berlinet04reproducing,saitoh16theory} allow embedding
probability measures into reproducing kernel Hilbert spaces (RKHS;
\citep{aronszajn50theory}) by use of a positive definite function, the
\textit{kernel function}. This approach has gained considerable attention over
the last 20 years. Such embeddings induce the so-called maximum mean
discrepancy (MMD; \citep{smola07hilbert,gretton12kernel}), which quantifies the
discrepancy of two probability measures by considering the RKHS norm of the
distance of their respective embeddings. MMD is
a metric on the space of probability distributions if the kernel is 
characteristic \citep{fukumizu08kernel,sriperumbudur10hilbert}. MMD is also an integral probability metric \citep{zolotarev83probability,muller97integral} where the underlying function class is chosen to be the unit ball in the corresponding RKHS.

MMD allows for the quantification of dependence by considering the distance between the embedding of a joint distribution and that of the product of its marginals. This construction gives rise to the so-called
Hilbert-Schmidt independence criterion (HSIC; \citep{gretton05measuring}), which
is also equal to the RKHS norm of the centered cross-covariance operator. In
fact, one of the most widely-used independence measures in statistics, distance
covariance \citep{szekely07measuring,szekely09brownian,lyons13distance}, was
shown to be equivalent to HSIC \citep{sejdinovic13equivalence} when the latter is specialized to $M=2$ components;
\citet{sheng23distance} proved a similar result for the conditional case. For
$M>2$ components
\citep{quadrianto09kernelized,sejdinovic13kernel,pfister18kernel}, universality
\citep{steinwart01influence,
micchelli06universal,carmeli10vector,sriperumbudur11universality} of the kernels
$(k_m)_{m=1}^M$ (on the respective domains) underlying HSIC guarantees that this
measure captures independence \citep{szabo18characteristic2}. In the case of
$M=2$, characteristic  $(k_m)_{m=1}^2$ suffice \citep{lyons13distance}.

HSIC has been deployed successfully in numerous
contexts, including independence testing in batch
\citep{gretton08kernel,wehbe15nonparametric,bilodeau17tests,goreczki18independence,pfister18kernel,albert22adaptive,shekhar23permutation}
and streaming \citep{podkopaev23sequential} settings, feature selection
\citep{camps10remote,song12feature,yamada14high,wang22rank} with applications in
biomarker detection \citep{gonzalez19block} and wind power prediction
\citep{bouche23wind}, clustering \citep{song07dependence,gonzalez19block}, and
causal discovery
\citep{mooij16distinguishing,pfister18kernel,chakraborty19distance,scholkopf21causal,kalinke23nystrom}.
In addition, HSIC has recently found successful applications in  sensitivity
analysis
\citep{veiga15global,freitas23sensitivity,fellmann24kernel,herrando24sensitivity},
in the context of uncertainty quantification \citep{stenger2020optimal}, for the
analysis of data augmentation methods for brain tumor detection
\citep{anaya22data}, and that of multimodal neural networks trained on
neuroimaging data \citep{fedorov24multimodal}.

Many
estimators for HSIC exist. The classical ones rely on U-statistics or
V-statistics \citep{gretton05measuring,quadrianto09kernelized,pfister18kernel}
and are known to converge at a rate of $\O_P\left(n^{-1/2}\right)$. In fact, the
V-statistic-based estimators are obtained by replacing the population kernel
mean embedding with its empirical counterpart; estimating the mean embedding can be carried out at a speed $\O_P\left(n^{-1/2}\right)$
\citep[Theorem~2]{smola07hilbert}, which implies that HSIC can be estimated at
the same rate. Existing approximations such as Nyström HSIC
\citep{kalinke23nystrom}, also achieve this rate under the assumption of an appropriate rate of decay of the effective dimension. While all of
these upper bounds match asymptotically, it is not known whether HSIC can be
estimated at a faster rate, that is, whether the upper bound of
$\O_P\left(n^{-1/2}\right)$ is optimal in the minimax sense, or if designing estimators
achieving better rates is possible. Lower bounds for the related MMD are known
\citep{tolstikhin16minimax}, but the existing analysis  considers radial kernels
and relies on independent Gaussian distributions.  Radial kernels are a special
case of the more general class of translation-invariant kernels that we
consider.\footnote{The family of radial kernels encompasses, for example, Gaussians, mixtures of
Gaussians, inverse multiquadratics, and Matérn kernels; the Laplace kernel is translation-invariant but not radial (with respect to the traditionally-chosen Euclidean norm $\left\|\cdot\right\|_{\R^d}$).}  The reliance on independent Gaussian distributions renders the
analysis of \citet{tolstikhin16minimax} inapplicable for HSIC estimation. We
tackle both of these severe restrictions in the present article.

We make the following \tb{contributions}.  
\begin{itemize}
\item  We establish the minimax lower bound $\O \left( n^{-1/2}
\right)$ of HSIC estimation with $M\ge 2$ components on $\R^{d}$ with continuous bounded translation-invariant characteristic kernels. As this
lower bound matches the known upper bounds of the existing ``classical''
U-statistic and V-statistic-based estimators, and that of the Nyström
HSIC estimator, our result settles their  minimax 
optimality.
\item Specifically, our result also implies the minimax lower bound of $\O \left(n^{-1/2}\right)$ for the estimation
of the cross-covariance operator, which can be further specialized to get back the minimax result \citep[Theorem~5]{zhou19covestimators} on the
estimation of the covariance operator.
\end{itemize} 

The paper is structured as follows. Notations are introduced in Section~\ref{sec:notations}. Section~\ref{sec:results} is dedicated to our main result on the minimax rate of HSIC estimation on $\R^d$, with proof presented in Section~\ref{sec:proofs}. An auxiliary result on the Kullback-Leibler divergence is shown in Appendix~\ref{sec:auxiliary-results}.
\section{Notations} \label{sec:notations}

In this section, we introduce a few notations $\mathbb N_{> 0}$, $[M]$, $\b
I_n$, $\bm 0_n$, $\bm 1_n$, $\b A\T$, $\ip{\b v,\b w}{}$, $\left\|\b v\right\|_{\R^d}$, $\bdiag\left(\b M_1,\ldots,\b M_N\right)$, $|\b A|$, $\mathcal
M_1^+\left(\R^d\right)$, $\psi_{\P}$, $\mathrm{KL}(\P||\Q)$, $L^2\left(\R^d,\Lambda\right)$, 
$\norm{f}{L^2\left(\R^d,\Lambda\right)}$, $\operatorname{supp}(\Lambda)$, $\H_k$,
$\phi_k$,  $k$, $\mu_k$, $\MMD_k$, $\otimes_{m=1}^M\H_{k_m}$, $\otimes_{m=1}^M k_m$,
$\P_m$, $\otimes_{m=1}^M\P_m$, $\P^n$, $\O_P\left(r_n\right)$, $\mathcal{O} (a_n)$, $a_n \asymp b_n$, $\HSIC_k$, and
$C_X$. Throughout the paper we consider random variables, probability measures, and kernels on~$\R^d$.

For $M \in \mathbb N_{>0} := \{1,2,\ldots\}$, let $[M] := \{1,\ldots,M\}$.
Denote by $\b I_n$ the $n\times n$-sized identity matrix and by $\bm 0_n =
(0,\ldots,0)\T\in\R^n$ (resp.\ $\bm 1_n = (1,\ldots,1)\T \in \R^n$) a column
vector of zeros (resp.\ ones). The transpose of a matrix $\b A \in \R^{d_1\times d_2}$ is written as $\b A\T \in \R^{d_2\times d_1}$. For $\b v, \b w\in\R^d$, $\ip{\b v,\b w}{} = \b v\T\b
w$ stands for their Euclidean inner product; $\left\|\b v\right\|_{\R^d}=\sqrt{\ip{\b v,\b v}{}}$ is the associated Euclidean norm. $\bdiag\left(\b M_1,\ldots,\b M_N\right)$ forms a block-diagonal
matrix from its arguments $(\b M_n)_{n=1}^N$ ($\b M_n\in \R^{d_n \times d_n}$, $n\in [N]$) and $|\b A|$ denotes the determinant of a matrix $\b A
\in \R^{d\times d}$.

The set of Borel probability measures on $\R^d$ is denoted by $\mathcal M_1^+ \left( \R^d \right)$.
For a random variable $X \sim \P \in \mathcal M_1^+\left(\R^d\right)$, we denote its
characteristic function by $\psi_{\P}(\bm\omega) = \E_{X\sim\P}\left[e^{i\ip{\bm
\omega,X}{}}\right]$ with $\bm \omega \in \R^d$ and $i=\sqrt{-1}$.  
Let $\P, \Q \in \mathcal M_1^+\left(\R^d\right)$, assume that $\P$ is absolutely continuous w.r.t.\ $\Q$, and let $\dfrac{\d \P}{\d \Q}$ denote the corresponding Radon-Nikodym derivative (of $\P$ w.r.t.\ $\Q$). Then, the Kullback-Leibler divergence of $\P$ and $\Q$ is defined as $\mathrm{KL}(\P||\Q) := \int_{\R^d}\log\left(\dfrac{\d \P}{\d \Q}(\b x)\right)\d \P(\b x)$.
Given a measure space $\left(\R^d,\mathcal B\left(\R^d\right), \Lambda\right)$, we denote by $L^2(\R^d,
\Lambda) := L^2\left(\R^d,\mathcal B\left(\R^d\right),\Lambda\right)$ the Hilbert space of
(equivalence classes of) measurable functions $f : \left(\R^d,\mathcal B\left(\R^d\right)\right) \to \left(\R,\mathcal B\left(\R\right)\right)$ for which
$\norm{f}{L^2\left(\R^d,\Lambda\right)}^{2} := \int_{\R^d}|f(\b x)|^2\d\Lambda(\b x) < \infty$. The
support of a probability measure $\Lambda \in \mathcal M_1^+\left(\R^d\right)$  denoted by $\operatorname{supp}(\Lambda)$ is the subset of
$\R^d$ for which every open neighborhood of $\b x\in\R^d$ has positive measure~\citep[p.~207]{cohn13measure}.

A function $k : \R^d \times \R^{d } \to \R$ is called a kernel if there exists a
Hilbert space $\H$ and a feature map $\phi : \R^{d }\to \H$ such that $k(\b x,
\b x') = \ip{\phi(\b x),\phi(\b x' )}{\H}$ for all $\b x, \b x' \in \R^d $. A Hilbert
space of functions $h :\R^d \to \R$ is an RKHS $\H_k$ associated to a kernel $k
:\R^d \times \R^d \to \R$ if $k(\cdot,\b x) \in \H_k$ and $\ip{h,k(\cdot,\b x)}{\H_k} =
h(\b x)$ for all $\b x \in \R^d$ and $h\in \H_k$.\footnote{For fixed $\b x\in\R^d$, the
function $k(\cdot,\b x) : \R^d \to\R$ means $\b x'\mapsto k(\b x',\b x)$.} 
In this work, we
assume all kernels to be measurable and bounded.\footnote{Boundedness of the
kernel, that is, $\sup_{\b x,\b x'\in\R^d}k(\b x,\b x') < \infty$, implies boundedness of the
feature map, that is, $\sup_{\b x\in\R^d}\norm{\phi_k(\b x)}{\H_k} < \infty$ (and vice
versa); it is also equivalent to $\sup_{\b x\in\R^d}k(\b x,\b x) < \infty$.} The function
$\phi_k(\b x) := k(\cdot,\b x)$ is the canonical feature map, and $k(\b x,\b x') =
\ip{k(\cdot,\b x),k(\cdot,\b x')}{\H_k} = \ip{\phi_k(\b x),\phi_k(\b x')}{\H_k}$ for all
$\b x, \b x'\in\R^d$. A function $\kappa :  \R^d \to \R$ is called positive definite if $\sum_{i,j\in[n]} c_ic_j\kappa(\b x_i -\b x_j) \ge 0$ for all $n\in\mathbb N_{>0}$, $\b c = \left(c_i\right)_{i=1}^n \in \R^n$, and $\{\b x_i\}_{i=1}^n \subset \R^d$. A kernel $k :\R^d\times \R^d\to \R$ is said to be
translation-invariant if there exists a positive definite function $\kappa : \R^d\to \R$ such that
$k(\b x,\b x') = \kappa(\b x-\b x')$ for all $\b x,\b x'\in \R^d$. By Bochner's theorem \citep[Theorem~6.6]{wendland05scattered} (recalled in Theorem~\ref{thm:bochner}) for a continuous bounded translation-invariant kernel $k:\R^d\times \R^d \to \R$ there exists a finite non-negative Borel measure $\Lambda_k$ such that
\begin{align}
k(\b x, \b y) &= \int_{\R^d} e^{-i \langle \b x-\b y, \bm \omega\rangle} \d \Lambda_k(\bm \omega) \label{eq:bochner}
\end{align}
for all $\b x, \b y \in \R^d$. The (kernel)
mean embedding of a probability measure $\P\in\mathcal M_1^+\left(\R^d\right)$ is
\begin{align*}
\mu_k(\P) = \int_{\R^d}\phi_k(\b x)\d\P(\b x) \in \H_k,
\end{align*}
where the integral is meant in Bochner's sense \citep[Chapter~II.2]{diestel77vector}; the boundedness of $k$ ensures that it is well-defined. For $\P,\Q \in \mathcal M_1^+\left(\R^d\right)$ one can define the (semi-)metric called maximum mean discrepancy \citep{smola07hilbert,gretton12kernel} as 
\begin{align*}
\MMD_k(\P,\Q) = \norm{\mu_k(\P)-\mu_k(\Q)}{\H_k}.
\end{align*}
If the mean embedding $\mu_k$ is injective, MMD is a metric and the kernel $k$
is called
characteristic~\citep{fukumizu08kernel,sriperumbudur10hilbert,szabo18characteristic2}.

Let $\R^d = \times_{m=1}^M\R^{d_m}$ ($d=\sum_{m=1}^Md_m$) and assume that each domain $\R^{d_m}$ is equipped with
a kernel $k_{m} : \R^{d_m}\times \R^{d_m}\to \R$ with associated RKHS
$\H_{k_{m}}$ ($m\in[M]$). The tensor product Hilbert space of
$\left(\H_{k_m}\right)_{m=1}^M$ is denoted by $\otimes_{m=1}^M\H_{k_m}$; it is
an RKHS \citep[Theorem~13]{berlinet04reproducing} with the tensor product kernel
$k = \otimes_{m=1}^Mk_m : \R^d \times \R^d \to \R$ defined by
\begin{align*}
k\left((\b x_m)_{m=1}^M, (\b x_m')_{m=1}^M\right) = \prod_{m\in[M]} k_m(\b x_m,\b x_m') \quad \text{for all}\quad \b x_m, \b x_m' \in \R^{d_m},\, m\in[M].
\end{align*}
The kernel $k$ has the canonical feature map $\phi_k\left((\b x_m)_{m=1}^M\right) = \otimes_{m=1}^M\phi_{k_m}\left(\b x_m\right) \in \otimes_{m=1}^M\H_{k_m} =: \H_k$ ($\b x_m \in \R^{d_m}, m\in[M]$).
Let $X = (X_m)_{m=1}^M$ be a random variable taking values in $\R^d$ with joint
distribution $\P \in \mathcal M_1^+\left(\R^d\right)$ and marginal distributions $\P_m \in
\mathcal M_1^+\left(\R^{d_m}\right)$ ($m\in[M]$; $d=\sum_{m=1}^M d_m$). We write $\otimes_{m=1}^M \P_m \in \mathcal
M_1^+\left(\R^d\right)$ for the product of measures $\P_m$ ($m\in [M]$). Specifically, $\P^n :=
\otimes_{i=1}^n\P \in \mathcal M_1^+\left(\left(\R^{d}\right)^n\right)$ denotes the $n$-fold product of $\P$. For a
sequence of real-valued random variables $\left(X_n\right)_{n=1}^\infty$ and a
sequence $\left(r_n\right)_{n=1}^\infty$ ($r_n>0$ for all $n$), $X_n = \O_P\left(r_n\right)$
denotes that $\frac{X_n}{r_n}$ is bounded in probability. For positive sequences $(a_n)_{n=1}^\infty$ and $(b_n)_{n=1}^\infty$, $b_n = \mathcal{O} (a_n)$ if there exist constants $C > 0$ and $n_0 \in \mathbb{N}_{>0}$ such that $b_n \le C a_n$ for all $n \ge n_0$; $a_n \asymp b_n$ if $a_n = \O\left(b_n\right)$ and $b_n = \O\left(a_n\right)$.  One can define our
quantity of interest, the Hilbert-Schmidt independence criterion (HSIC;
\citep{gretton05measuring,quadrianto09kernelized,pfister18kernel,szabo18characteristic2}), as
\begin{align}
  \HSIC_k(\P) &= \MMD_k\left(\P,\otimes_{m=1}^M\P_m\right) = \norm{C_X}{\H_k},\\
  C_X &= \mu_k(\P)-\mu_k\left(\otimes_{m=1}^M\P_m\right) \in \H_k,\label{eq:c-x}
\end{align}
and $C_{X}$
denotes the centered cross-covariance operator.

\section{Results} \label{sec:results}

This section is dedicated to our results: The minimax lower bound for the
estimation of $\HSIC_k(\P)$, where $k$ is a product of continuous bounded translation-invariant characteristic kernels
is given in Theorem~\ref{thm:main-stmt}(ii). For the specific case where $k$ is a product of Gaussian kernels (stated in Theorem~\ref{thm:main-stmt}(i)), the constant in the lower bound is made explicit.  Theorem~\ref{thm:main-stmt}(ii) also helps to establish a lower bound on the estimation of the cross-covariance operator (Corollary~\ref{corollary:covariance-estimation}). 

Before presenting our results, we recall the framework of minimax estimation \citep{tsybakov09nonparametric} adapted to our setting. Let $\hat F_n$ denote any estimator of $\HSIC_k(\P)$ based on $n$ i.i.d.\ samples from $\P$. A sequence $(\xi_n)_{n=1}^\infty$ ($\xi_{n} > 0$ for all $n$) is said to be a lower bound of HSIC estimation w.r.t.\  a class $\mathcal P$ of Borel probability measures on $\R^d$ if there exists a constant $c>0$ such that
\begin{align}
  \inf_{\hat F_n}\sup_{\P\in\mathcal P}\P^n\left\{\xi_n^{-1}\left|\HSIC_k(\P)-\hat F_n\right|\ge c \right\}  > 0. \label{eq:prob-bounds}
\end{align}
If a specific estimator of HSIC $\tilde F_n$ has an upper bound that matches $\left( \xi_{n} \right)_{n=1}^{\infty}$ up to constants, that is,
\begin{align} \label{eq:upper-bounds}
  \left| \HSIC_{k}(\P) - \tilde F_{n} \right| = \O_{P}\left(\xi_{n} \right),
\end{align}
then $\tilde F_n$ is called minimax optimal.

We use Le Cam's method \citep{lecam73convergence,tsybakov09nonparametric}
(recalled in Theorem~\ref{theorem:le-cam}) to obtain bounds as in
\eqref{eq:prob-bounds}; estimators of HSIC achieving the bounds in
\eqref{eq:upper-bounds} with $\xi_{n} = n^{-1/2}$ are quoted in the
introduction. 
The key to the application of the method is to show that 
there exist $\alpha>0$ and $n_0 \in \mathbb{N}_{>0}$ such that for all $n\ge n_0$ one can find an adversarial pair of
distributions $\left(\P_{\theta_0},\P_{\theta_1}\right) = \left(\P_{\theta_0}(n),\P_{\theta_1}(n)\right) \in \mathcal P \times \mathcal P$ and $s_n >0$  for which 
\begin{enumerate}
\item $\operatorname{KL}\left(\P_{\theta_1}^n||\P_{\theta_0}^n\right) \le \alpha$, in other words, the corresponding $n$-fold product measures must be similar in the sense of Kullback-Leibler divergence, but 
\item  $\left|\HSIC_{k}(\P_{\theta_1})-\HSIC_k(\P_{\theta_0})\right| \ge 2s_n$, that is, their corresponding values of HSIC must be dissimilar. 
\end{enumerate}
In this case, $\inf_{\hat F_n} \sup_{\P\in\mathcal P} \P^n\left\{\left|\HSIC_k(\P)-\hat F_n\right|\ge s_n \right\}  \ge \max \left(\frac{e^{-\alpha}}{4},\frac{1-\sqrt{\alpha/2}}{2}\right)$ for all $n\ge n_0$; hence to establish the minimax optimality of existing estimators w.r.t.\
their known upper bounds, it is sufficient to find adversarial pairs $\left\{\left(\P_{\theta_0}(n),\P_{\theta_1}(n)\right)\right\}_{n\ge n_0}$ that satisfy 1.\ for some positive constant $\alpha$ and also fulfill 2.\ with $s_n\asymp n^{-1/2}$.

The proof of the first part of our statement relies on the following Lemma~\ref{lemma:closed-form-hsic} which yields the analytical value of
$\HSIC_k\left(\mathcal N(\bm \mu,\bm \Sigma)\right)$,  where $k = \otimes_{m=1}^M
k_m$ is the product of Gaussian kernels $k_m$ ($m\in[M]$) and $\mathcal N(\bm
\mu,\bm \Sigma)$ denotes the multivariate normal distribution with mean $\bm\mu
\in \R^d$ and covariance matrix $\bm\Sigma\in\R^{d\times d}$.

\begin{lemma}[Analytical value of HSIC for the Gaussian setting]
  \label{lemma:closed-form-hsic}
Let us consider the Gaussian kernel $k(\b x,\b y) = e^{-\frac{\gamma}{2} \left \| \b x - \b y \right\|_{\R^d}^2}$ ($\gamma>0$, $\b x, \b y \in \R^d$) and Gaussian random variable $X = (X_m)_{m=1}^M \sim \mathcal N(\b m, \bm \Sigma) =: \P$, where $X_m \in \R^{d_m}$ ($m\in[M]$), $\b m = (\b m_m)_{m=1}^M \in \R^{d}$, $\bm \Sigma = [\bm \Sigma_{i,j}]_{i,j\in [M]} \in \R^{d \times d}$, $\bm \Sigma_{i,j} \in \R^{d_i \times d_j}$, and $d=\sum_{m\in [M]}d_m$.  
In this case, with $\bm \Sigma_1 = \bm\Sigma$ and $\bm\Sigma_2 = \bdiag(\bm \Sigma_{1,1},\dots,\bm \Sigma_{M,M})$, we have
\begin{align*}
    \HSIC^2_k(\P) &=  \frac{1}{\left|2\gamma\bm{\Sigma}_1+\b{I}_d\right|^{\frac{1}{2}}} +  \frac{1}{\left|2\gamma\bm{\Sigma}_2+\b{I}_d\right|^{\frac{1}{2}}} - \frac{2}{\left|\gamma\bm{\Sigma}_1+\gamma\bm{\Sigma}_2+\b{I}_d\right|^{\frac{1}{2}}}.
\end{align*}
\end{lemma}

In this work, we focus on continuous bounded translation-invariant kernels, which are
fully characterized by Bochner's theorem
\citep[Theorem~6.6]{wendland05scattered}; the theorem states that a function on $\R^d$
is positive definite if and only if it is the Fourier transform of a finite
nonnegative measure.\footnote{We note that for
many translation-invariant kernels, the corresponding spectral measures are
known \citep[Table~2]{sriperumbudur10hilbert}.} We use this description to obtain our main result, which is as follows.

\begin{theorem}[Lower bound for HSIC estimation on $\R^d$]
  \label{thm:main-stmt}
  Let $\mathcal P$ be a class of Borel probability measures over $\R^d$
  containing the $d$-dimensional Gaussian distributions. Let
  $d=\sum_{m\in[M]}d_m$ and $\hat F_n$ denote any estimator of $\HSIC_k(\P)$
  with $n \ge 2 =: n_0$ i.i.d.\ samples from $\P \in \mathcal P$. Assume further that
  $k=\otimes_{m=1}^M k_m$ where either, for $m\in[M]$,
  \begin{enumerate}[label=(\roman*)]
    \item the kernels $k_m : \R^{d_m}\times \R^{d_m} \to \R$ are Gaussian with common bandwidth parameter $\gamma
    >0$ defined by $\left(\b x_m, \b x_m'\right)
    \mapsto e^{-\frac{\gamma}{2}\norm{\b x_m-\b x_m'}{\R_{d_m}}^2}$ ($\b x_m, \b x_m' \in \R^{d_m}$),  or
    \item the kernels $k_m : \R^{d_m}\times \R^{d_m} \to \R$ are continuous bounded     translation-invariant characteristic kernels.
  \end{enumerate}
  Then, for any $n\ge n_0$, it holds that
    \begin{align*}
      \inf_{\hat F_n}\sup_{\P\in\mathcal P}\P^n\left\{\left|\HSIC_k\left(\P\right)-\hat F_n\right|\ge \frac{c}{\sqrt n}\right\} \ge \frac{1-\sqrt{\frac{5}{8}}}{2},
    \end{align*}
    with  (i) the constant $c = \frac{\gamma}{ 2\left(2\gamma+1\right)^{\frac{d}{4}+1}}>0$ (depending on $\gamma$ and $d$ only) in the first case, or
      (ii) some constant $c > 0$ in the second case.
  \end{theorem}

We note that while Theorem~\ref{thm:main-stmt}(ii) applies to the more general
class of translation-invariant kernels, we include Theorem~\ref{thm:main-stmt}(i) as it makes the constant $c$ explicit.

The following corollary allows to recover the recent lower bound on the
estimation of the covariance operator by \citet[Theorem~5]{zhou19covestimators}
as a special case that we detail in Remark~\ref{remark:main}(e).

\begin{corollary}[Lower bound on cross-covariance operator estimation] \label{corollary:covariance-estimation}
  In the setting of Theorem~\ref{thm:main-stmt}(ii), let $\hat F_n$ denote any estimator of the centered cross-covariance operator $C_X \in \H_{k}$ defined in~\eqref{eq:c-x} with $n\ge 2 =: n_0$ i.i.d.\ samples from $\P \in \mathcal P$. Then, for any $n\ge n_0$, it holds that
  \begin{align*}
    \inf_{\hat F_n}\sup_{\P\in\mathcal P}\P^n\left\{\norm{C_X - \hat F_n}{\H_k}\ge \frac{c}{\sqrt n}\right\} \ge \frac{1-\sqrt{\frac{5}{8}}}{2},
  \end{align*}
  for some constant $c >0 $.
\end{corollary}

\begin{remark}~
  \label{remark:main}
  \begin{enumerate}[label=(\alph*)]
  \item \tb{Validness of HSIC.} Though generally the characteristic property of $(k_m)_{m=1}^M$-s is not enough  \citep[Example~2]{szabo18characteristic2} for $M> 2$ to ensure the $\mathcal{I}$-characteristic property of $k=\otimes_{m=1}^M k_m$ (in other words, that $\text{HSIC}_{k}(\P)=0$ iff.\ $\P=\otimes_{m=1}^M \P_m$), on $\R^d$ under the imposed continuous bounded translation-invariant assumption (i) $k$ being characteristic, (ii) $k$ being $\mathcal{I}$-characteristic, and (iii) $(k_m)_{m=1}^M$-s being characteristic are equivalent (Theorem~\ref{thm:szabo-thm4}).
  \item \tb{Minimax optimality of existing HSIC estimators.} The lower bounds in
   Theorem~\ref{thm:main-stmt} asymptotically match the known upper bounds of
   the U-statistic and V-statistic-based estimators of $\xi_{n} = n^{-1/2}$. The
   Nyström-based HSIC estimator achieves the same rate under an appropriate
   decay of the eigenspectrum of the respective covariance operator.  Hence, Theorem~\ref{thm:main-stmt}
   implies the optimality of these estimators on $\R^{d}$
   with continuous bounded translation-invariant characteristic kernels in the minimax sense.

  \item \tb{Difference compared to \citet{tolstikhin16minimax} (minimax MMD estimation).}  We note that a
    lower bound for the related MMD$_k$ exists. However, the adversarial
    distribution pair $(\P_{\theta_1}, \P_{\theta_0})$ constructed by
    \citet[Theorem~1]{tolstikhin16minimax} to obtain the lower bound on MMD
    estimation has a product structure which implies that
    $\left|\HSIC_{k}(\P_{\theta_1})-\HSIC_k(\P_{\theta_0})\right| = 0$
    and hence it is not applicable in our case of HSIC;
    \citet[Theorem~2]{tolstikhin16minimax} with radial kernels has the same
    restriction.

  \item \tb{Difference compared to \citet{tolstikhin17minimax} (minimax mean embedding estimation).} The estimation
  of the mean embedding  $\mu_k(\P)$ is known to have a minimax rate of $\O
  \left(n^{-1/2}\right)$. But, this rate does not imply an optimal lower bound for the
  estimation of MMD as is evident from the two works \citep{tolstikhin16minimax,tolstikhin17minimax}. The same conclusion holds for HSIC
  estimation. 

  \item \tb{Difference compared to \citet{zhou19covestimators} (minimax covariance operator estimation).}
    For the related problem of estimating the centered \textit{covariance} operator 
    \begin{align*}
      C_{\mathit{XX}} = \int_{\R^d}\left(\phi_k(x)-\mu_k(\P)\right)\otimes\left(\phi_k(x)-\mu_k(\P)\right)\d\P(x) \in \H_k \otimes \H_k,
    \end{align*}
    \citet[Theorem~5]{zhou19covestimators} give the lower bound
    \begin{align*}
      \inf_{\hat F_n}\sup_{\P\in\mathcal P}\P^n\left\{\norm{C_{\mathit{XX}} -\hat F_n}{\H_k\otimes\H_k} \ge \frac{c}{\sqrt{n}}\right\} \ge 1/8
    \end{align*}
    in the same setting as in Theorem~\ref{thm:main-stmt}(ii), where $\hat F_n$
    is any estimator of the centered covariance $C_{\mathit{XX}}$, and $c$ is
    a positive constant. By noting that the centered covariance is the centered
    cross-covariance of a random variable with itself,
    Corollary~\ref{corollary:covariance-estimation} recovers their result.
  \end{enumerate}
\end{remark}
The next section contains our proofs.

\section{Proofs} \label{sec:proofs}

This section is dedicated to our proofs. We present the proof of Lemma~\ref{lemma:closed-form-hsic} in Section~\ref{proof:analytic-hsic},
that of Theorem~\ref{thm:main-stmt} in Section~\ref{sec:proof-main-stmt}, and that of
Corollary~\ref{corollary:covariance-estimation} in
Section~\ref{sec:proof-corollary}. 

\subsection{Proof of Lemma~\ref{lemma:closed-form-hsic}}\label{proof:analytic-hsic}
    As 
    \begin{align*}
    \HSIC^2_k(\P)  &=\MMD^2_k(\P,\Q) = \left\|\mu_k(\P) - \mu_k(\Q)\right\|_{\H_k}^2 \\
    & = \langle \mu_k(\P), \mu_k(\P)\rangle_{\H_k} + \langle \mu_k(\Q), \mu_k(\Q)\rangle_{\H_k} - 2 \langle \mu_k(\P),\mu_k(\Q) \rangle_{\H_k}
    \end{align*}
    with $\Q = \otimes_{m=1}^M \P_m= \mathcal N(\b m, \bdiag(\bm
     \Sigma_{1,1},\dots,\bm \Sigma_{M,M}))$, $\P_m = \mathcal N(\b m_m,\bm
     \Sigma_{m,m})$, it is sufficient to be able to compute $\langle
     \mu_k(\P),\mu_k(\Q)\rangle_{\H_k}$-type quantities with $\P = \mathcal N(\b
     m_1, \bm \Sigma_1)$ and $\Q = \mathcal N(\b m_2, \bm \Sigma_2)$. One can
     show \citep[Table~1]{muandet12learning} that
$\langle \mu_k(\P),\mu_k(\Q)\rangle_{\H_k}  = \frac{e^{-\frac{1}{2}\left(\b{m}_1-\b{m}_2\right)\T \left(\bm{\Sigma}_1+\bm{\Sigma}_2+\gamma^{-1}\b{I}_d\right)^{-1}\left(\b{m}_1-\b{m}_2\right)}}{\left|\gamma\bm{\Sigma}_1+\gamma\bm{\Sigma}_2+\b{I}_d\right|^{\frac{1}{2}}}$.
Using this fact and that $\b m = \b m_1 = \b m_2$, the result follows.

\subsection{Proof of Theorem~\ref{thm:main-stmt}}
\label{sec:proof-main-stmt}

The setup and the upper bound on
$\operatorname{KL}(\P_{\theta_1}^n||\P_{\theta_0}^n)$ agree for (i) and (ii) but
the methods that we use to lower bound
$\left|\HSIC_{k}(\P_{\theta_1})-\HSIC_k(\P_{\theta_0})\right|$ differ. We
structure the proof accordingly and present the overlapping part before we
branch out into (i) and (ii). Both parts of the statement rely on Le Cam's
method, which we state as Theorem~\ref{theorem:le-cam} for self-completeness. 

To construct the adversarial pair, we consider a class $\mathcal G$ of Gaussian
distributions over $\R^d$ such that every element $\mathcal{N}\big(\bm \mu, \bm
\Sigma \big) \in \mathcal G$, with
\begin{align}
\bm \Sigma = \bm \Sigma(i,j,\rho) =
       \begin{bmatrix}   1   & \cdots & 0      & 0      & \cdots & 0      \\
                      \vdots & \ddots & \vdots & \vdots &        & \vdots \\
                         0   & \cdots & 1      & \rho   & \cdots & 0      \\
                         0   & \cdots & \rho   & 1      & \cdots & 0      \\
                      \vdots &        & \vdots & \vdots & \ddots & \vdots \\
                         0   & \cdots & 0      & 0      & \cdots & 1
       \end{bmatrix} \in \R^{d\times d}, \label{eq:Sigma-rho}
\end{align}
and (fixed) $i = d_1$, $j=d_1+1$, $\rho\in(-1,1)$. In other words, $\bm \Sigma$
is essentially the $d$-dimensional matrix $\bm I_d$ except for the $(i,j)$ and
$(j,i)$ entry; both entries are identical to $\rho$, and they specify the correlation of
the respective coordinates. This family of distributions is indexed by a tuple
$(\bm \mu,\rho) \in \R^d  \times (-1,1) =:  \mathcal A$ and, for $a \in \mathcal A$, we write $\P_a$ for the associated distribution.  To bring
ourselves into the setting of Theorem~\ref{theorem:le-cam}, we fix $n\in\mathbb N_{>0}$, choose $\X =  \left(\R^d\right)^n$, set $\Theta  = \{\theta_{a} := \HSIC_k(\P_{a}) \,:\, a \in \mathcal A \}$, $\mathcal P_{\Theta} = \{\P_a^n \,:\, a \in \mathcal A\} = \{\P_a^n \,:\, \theta_a \in \Theta\}$, and use the metric $(x,y) \mapsto |x-y|$ for
$x,y\in\R$. Hence, the data $D\sim
\P_\theta \in \mathcal P_\Theta$. For brevity, let $F : \mathcal A \to \R$ stand for $a \mapsto \HSIC_k(\P_a)$, and let $\hat F_n$ stand for the corresponding estimator based on $n$ samples.

As $\mathcal G \subseteq \mathcal P$, it holds for every positive $s$ that
  \begin{align*}
    \sup_{\P\in\mathcal P}\P^n\left\{\left|\HSIC_k\left(\P\right)-\hat F_n\right|\ge s \right\} \ge \sup_{\P\in\mathcal G}\P^n\left\{\left|\HSIC_k\left(\P\right)-\hat F_n\right|\ge s \right\}.
  \end{align*}

  Let $\P_{\theta_0} = \mathcal{N}\left(\bm \mu_0, \bm \Sigma_0\right)$ and $\P_{\theta_1} = \mathcal{N}\left(\bm \mu_1, \bm \Sigma_1\right)$ with 
  \begin{align*}
    \bm\mu_0 &= \bm0_d \in \R^d,
    & \bm \Sigma_0 &= \bm  \Sigma(d_1,d_1+1,0) = \mathbf I_d \in \R^{d\times d}, \\
    \bm\mu_1 &= \frac{1}{\sqrt dn} \bm 1_d \in \R^d,
    & \bm \Sigma_1 &= \bm \Sigma(d_1,d_1+1,\rho_n) \in \R^{d\times d},
  \end{align*}
  where $\rho_n \in (-1,1)$ will be chosen appropriately later.\footnote{Notice the dependence of $\P_{\theta_1}$ on $n$.}  We now proceed to upper bound
  $\mathrm{KL}\left(\P_{\theta_1}^n||\P_{\theta_0}^n\right)$ and lower bound
  $|F(\theta_1)-F(\theta_0)|$.

\paragraph{Upper bound for KL divergence} 
Lemma~\ref{lemma:bound-kl} implies that with $\rho_n^2=\frac{1}{n}$, one has the bound $\mathrm{KL}\left(\P_{\theta_1}^n||\P_{\theta_0}^n\right) \le \alpha:=\frac{5}{4}$ for $n\ge 2 =: n_0$. 

\paragraph{Lower bound (i): Gaussian kernels.}
Recall that the considered kernel is $k(\b x, \b y)=e^{-\frac{\gamma}{2} \left\|\b x - \b y \right\|_{\R^d}^2}$ ($\gamma >0$). 
The idea of the proof is as follows.
\begin{enumerate}
    \item We express $|F(\theta_1)-F(\theta_0)|$ in closed form as a function of $\gamma$, $\rho_n$, and $d$.
    \item Using the analytical form obtained in the 1st step, we construct the lower bound. 
\end{enumerate}
This is what we detail next.

\begin{itemize}
    \item \tb{Analytical form of $|F(\theta_1)-F(\theta_0)|$}:
    Using the fact that $\HSIC_k(\P_{\theta_0}) = 0$, we have that
\begin{eqnarray*}
  \lefteqn{\big|F\left(\theta_1\right)-\underbrace{F\left(\theta_0\right)}_{=0}\big|^2 = F^2\left(\theta_1\right)= \HSIC_k^2\left(\P_{\theta_1}\right)  = \MMD_k^2\left(\mathcal N\left(\bm\mu_1,\bm\Sigma_1\right), \mathcal N\left(\bm\mu_1,\b I_d\right)\right)} \\
  && \hspace*{-0.65cm} =  \norm{\mu_k\left(\mathcal N\left(\bm\mu_1,\bm\Sigma_1\right)\right) - \mu_k\left(\mathcal N\left(\bm\mu_1,\b I_d\right)\right)}{\H_k}^2 \\
  && \hspace*{-0.65cm} = \underbrace{\ip{\mu_k\left(\mathcal N\left(\bm\mu_1,\bm\Sigma_1\right)\right),\mu_k\left(\mathcal N\left(\bm\mu_1,\bm\Sigma_1\right)\right)}{\H_k}}_{(i)} + \underbrace{\ip{\mu_k\left(\mathcal N\left(\bm\mu_1,\b I_d\right)\right),\mu_k\left(\mathcal N\left(\bm\mu_1,\b I_d\right)\right)}{\H_k}}_{(ii)} \\
  && \hspace*{-0.65cm} \quad- 2\underbrace{\ip{\mu_k\left(\mathcal N\left(\bm\mu_1,\bm\Sigma_1\right)\right),\mu_k\left(\mathcal N\left(\bm\mu_1,\b I_d\right)\right)}{\H_k}}_{(iii)},
\end{eqnarray*}
which we compute term-by-term with Lemma~\ref{lemma:closed-form-hsic}, and obtain
\begin{align*}
  (i) &= \left|2 \gamma\bm \Sigma_1 + \b I_d\right|^{-1/2} = \left[\left(2\gamma+1\right)^{d-2}\left(\left(2\gamma+1\right)^2-\left(2\gamma\rho_n\right)^2\right)\right]^{-1/2},\\
  (ii) &= \left|2\gamma \b I_d+\b I_d\right|^{-1/2} = \left[\left(2\gamma +1\right)^d\right]^{-1/2}, \\
  (iii) &= \left|\gamma\bm\Sigma_1 + \gamma \b I_d + \b I_d\right|^{-1/2} = \left[\left(2\gamma+1\right)^{d-2}\left(\left(2\gamma+1\right)^2-\left(\gamma\rho_n\right)^2\right)\right]^{-1/2}.
\end{align*}
Combining (i), (ii), and (iii) yields that
\begin{align*}
  \HSIC_k^2\left(\P_{\theta_1}\right) &= (i) + (ii)-2(iii)  \\
  &= \left[\left(2\gamma+1\right)^{d-2}\left(\left(2\gamma+1\right)^2-\left(2\gamma\rho_n\right)^2\right)\right]^{-1/2} +
  \left[\left(2\gamma +1\right)^d\right]^{-1/2}\\
  &\quad -2 \left[\left(2\gamma+1\right)^{d-2}\left(\left(2\gamma+1\right)^2-\left(\gamma\rho_n\right)^2\right)\right]^{-1/2}.
\end{align*}
    \item \tb{Lower bound on $|F(\theta_1)-F(\theta_0)|$}:
    Next, we show that there exists $c>0$ such that for any $n \in \mathbb N_{>0}$ it holds that
$\HSIC_k^2\left(\P_{\theta_1}\right) \ge \frac cn$.

For $0< x < \left(1+\frac{1}{2\gamma}\right)^2$, let us consider the function
\begin{align*}
f_{c}(x) &=
       \left[\left(2\gamma+1\right)^{d-2}\left(\left(2\gamma+1\right)^2-4\gamma^2x\right)\right]^{-1/2}+\left[\left(2\gamma+1\right)^d\right]^{-1/2}\\
     &\quad-2\left[\left(2\gamma+1\right)^{d-2}\left(\left(2\gamma+1\right)^2-\gamma^2x\right)\right]^{-1/2}-cx \\
     &= \left[z^{d-2}\left(z^2-4\gamma^2x\right)\right]^{-1/2}+\left(z^d\right)^{-1/2}-2\left[z^{d-2}\left(z^2-\gamma^2x\right)\right]^{-1/2}-cx,
\end{align*}
with the shorthand $z:=2\gamma+1$.\footnote{Notice that $\left(2\gamma+1\right)^2-\gamma^2x>\left(2\gamma+1\right)^2-4\gamma^2x$, and $\left(2\gamma+1\right)^2-4\gamma^2x>0 \Leftrightarrow x<\left(1+\frac{1}{2\gamma}\right)^2$ for a positive $x$; hence the imposed assumption on $x$ ensures that the function $f_c$ is well-defined.}
With this notation, $f_{c}(1/n) = \HSIC_k^2\left(\P_{\theta_1}\right) -c/n$; our aim is to determine $c>0$ such that $f_{c}(1/n) \ge 0$ for any positive integer $n$. To achieve this goal, notice that $f_{c}(0)=0$, and
\begin{align*}
  f_{c}'(x)&=  \frac{2\gamma^2z^{d-2}}{\left[z^{d-2}\left(z^2-4x\gamma^2\right)\right]^{3/2}}-\frac{\gamma^2z^{d-2}}{\left[z^{d-2}\left(z^2-x\gamma^2\right)\right]^{3/2}}-c \\
       &> \frac{2\gamma^2z^{d-2}}{\left[z^{d-2}\left(z^2-x\gamma^2\right)\right]^{3/2}}-\frac{\gamma^2z^{d-2}}{\left[z^{d-2}\left(z^2-x\gamma^2\right)\right]^{3/2}}-c = \frac{\gamma^2z^{d-2}}{\left[z^{d-2}\left(z^2-x\gamma^2\right)\right]^{3/2}}-c \\
  &> \frac{\gamma^2z^{d-2}}{\left(z^{d-2}z^2\right)^{3/2}}-c = \frac{\gamma^2}{z^2\sqrt{z^d}}-c = \frac{\gamma^2}{\left(2\gamma+1\right)^2\sqrt{\left(2\gamma+1\right)^d}} -c.
\end{align*}
Choosing now $c=\frac{\gamma^2}{\left(2\gamma+1\right)^2\sqrt{\left(2\gamma+1\right)^d}}>0$, we have $f_{c}'(x) \ge 0$, so $f$ is a nondecreasing function. Note that $f_{c}(1/n) = \HSIC_k^2\left(\P_{\theta_1}\right) -c/n \ge 0$, with $x=1/n$ and ${\left(1+\frac{1}{2\gamma}\right)^{-2}} <1 \le n < \infty$.
By taking the positive square root, this means that 
\begin{align}
    \HSIC_k\left(\P_{\theta_1}\right) \ge \frac{\gamma}{ \left(2\gamma+1\right)\left(\left(2\gamma+1\right)^d\right)^{1/4}\sqrt n}=:2s
\end{align}
holds for $n\ge 1$,
implying that $|F(\theta_1)-F(\theta_0)| \ge 2s > 0$.
\end{itemize}

We conclude the proof by Theorem~\ref{theorem:le-cam} using that $\alpha=\frac{5}{4}$ and 
$\max\left(\frac{e^{-\frac{5}{4}}}{4}, \frac{1-\sqrt{\frac{5}{8}}}{2} \right) = \frac{1-\sqrt{\frac{5}{8}}}{2}$.

\paragraph{Lower bound (ii): translation-invariant kernels.}
Let $\Lambda_k$ denote the spectral measure associated to the kernel $k$ according to \eqref{eq:bochner}.
Using the fact that $\HSIC_k(\P_{\theta_0}) = 0$, we have for $|F(\theta_1)-F(\theta_0)|$ that
\begin{eqnarray}
  \lefteqn{\big|F\left(\theta_1\right)-\underbrace{F\left(\theta_0\right)}_{=0}\big|^2
  = F^2\left(\theta_1\right)= \HSIC_k^2\left(\P_{\theta_1}\right)  =\MMD_k^2\left(\mathcal N\left(\bm\mu_1,\bm\Sigma_1\right), \mathcal N\left(\bm\mu_1,\bm\Sigma_0\right)\right)}
  \\
  &&\stackrel{(i)}{=} \norm{\psi_{\mathcal N(\bm\mu_1,\bm \Sigma_1)}-\psi_{\mathcal N(\bm\mu_1,\bm\Sigma_0)}}{L^2\left(\R^d,\Lambda_k\right)}^2 \\
  &&\stackrel{(ii)}{=}  \int_{\R^d}  \left|e^{i\ip{\bm\mu_1,\bm \omega}{}-\frac12\ip{\bm\omega,\bm\Sigma_1\bm\omega}{}}-e^{i\ip{\bm\mu_1,\bm \omega}{}-\frac12\ip{\bm\omega,\bm\Sigma_0\bm\omega}{}}\right|^2\d\Lambda_k(\bm\omega) \\
  &&\stackrel{\hphantom{(ii)}}{=} \int_{\R^d}  \underbrace{\left|e^{i\ip{\bm\mu_1,\bm \omega}{}}\right|^2}_{=1}\left|e^{-\frac12\ip{\bm\omega,\bm\Sigma_1\bm\omega}{}}-e^{-\frac12\ip{\bm\omega,\bm\Sigma_0\bm\omega}{}}\right|^2\d\Lambda_k(\bm\omega) \\
  && \stackrel{(iii)}{\ge} \int_A  \left|e^{-\frac12\ip{\bm\omega,\bm \Sigma_1\bm\omega}{}}-e^{-\frac12\ip{\bm\omega,\bm \Sigma_0\bm\omega}{}}\right|^2\d\Lambda_k(\bm\omega) 
  \stackrel{(iv)}{\ge} \rho_n^2\underbrace{\int_A  \left[h'_{\bm\omega}(0)\right]^2\d\Lambda_k(\bm\omega)}_{=: (2c)^2} 
  \stackrel{(v)}{=} \underbrace{\frac{(2c)^2}{n}}_{=: (2s)^2 > 0}, 
\end{eqnarray}
where
$(i)$ holds by \citet[Corollary~4(i)]{sriperumbudur10hilbert} (recalled in Theorem~\ref{thm:bharath-char-func}).
$(ii)$ follows from the analytical form $\psi_{\mathcal N(\bm \mu, \bm \Sigma)}(\b t) = e^{i\ip{\bm \mu, \bm t}{} -\frac{1}{2}\ip{\b t,\bm\Sigma \b t}{}}$ of the  characteristic function of a multivariate normal distribution $\mathcal N(\bm \mu,\bm \Sigma)$.
For $(iii)$, we define the non-empty open set
\begin{align*}
A = \left\{ \bm\omega = (\omega_1,\ldots,\omega_d)\T \in \R^d \,:\, \omega_{d_1}\omega_{d_1+1} < 0\right\} \subset \R^d,
\end{align*}
and use that the integration of a non-negative function over a subset yields a lower bound.
 In $(iv)$, fix $\bm\omega \in A$ and let 
 \begin{align*}
 h_{\bm \omega} : \rho \in [0,1] \mapsto  e^{-\frac{1}{2}\ip{\bm\omega,\bm\Sigma(d_1,d_1+1,\rho)\bm\omega}{}} \in (0,1].
 \end{align*}
 Note that $h_{\bm\omega}(\rho) = e^{-\frac{1}{2}\left(\bm\omega\T\bm\omega + 2\rho\omega_{d_1}\omega_{d_1+1}\right)}$; $h_{\bm\omega}$ is continuous on $[0,1]$ and differentiable on $(0,1)$. Hence for any $\rho \in (0,1)$, by the mean value theorem, there exists $\tilde \rho\in(0,1)$ such that
  \begin{align*}
    h_{\bm \omega}(\rho) - h_{\bm \omega}(0) = \rho h_{\bm \omega}'( \tilde \rho) \ge \rho \min_{c\in[0,1]}h_{\bm \omega}'(c).
  \end{align*}
  We have the first and second derivatives 
  \begin{align*}
  h_{\bm\omega}'(c) &=
  -\omega_{d_1}\omega_{d_1+1}e^{-\frac{1}{2}\left(\bm\omega\T\bm\omega +
  2c\omega_{d_1}\omega_{d_1+1}\right)},& 
  h_{\bm\omega}''(c) &=
  \omega_{d_1}^2\omega_{d_1+1}^2e^{-\frac{1}{2}\left(\bm\omega\T\bm\omega +
  2c\omega_{d_1}\omega_{d_1+1}\right)} >0,
  \end{align*}
  which implies that  $c \mapsto
  h_{\bm \omega}'(c)$ is a strictly increasing function of $c$ and that it
  attains its minimum at $c = 0$, that is,
  \begin{align*}
  h_{\bm \omega}(\rho) - h_{\bm
  \omega}(0) \ge \rho h_{\bm \omega}'(0) >0,
  \end{align*}
  where the 2nd inequality holds by $\rho>0$ and $\bm \omega \in A$.
  This shows that 
  \begin{align}
      \left[h_{\bm \omega}(\rho) - h_{\bm
  \omega}(0)\right]^2 \ge \left[\rho h_{\bm \omega}'(0)\right]^2,
  \end{align}
  and the monotonicity of integration gives $(iv)$.
  For  $(v)$, we note that
   the kernel $k=\otimes_{m=1}^Mk_m$ is characteristic \citep[Theorem~4]{szabo18characteristic2} (recalled in
  Theorem~\ref{thm:szabo-thm4}) as the $(k_m)_{m=1}^M$-s are characteristic. Thus,
  $\operatorname{supp}\left(\Lambda_k\right) = \R^d$ (see
  \citet[Theorem~9]{sriperumbudur10hilbert}; recalled in
  Theorem~\ref{thm:bharath-full-Rd}), implying that
  $\Lambda_k(A) > 0$. $(v)$ follows from the
  positivity of $h_{\bm\omega}'(0)$ (for any $\bm\omega\in A$), from the fact that
  the integral of a positive function on a set with positive measure is
  positive, and from our choice of $\rho_n = n^{-1/2}$.

Now, by taking the positive square root, we have
\begin{align}
    \left|F\left(\theta_1\right)-F\left(\theta_0\right)\right| \ge \frac{2c}{\sqrt
n} =: 2s. \label{eq:ti-bound}
\end{align}
We conclude by the application of Theorem~\ref{theorem:le-cam}
using that $\alpha=\frac{5}{4}$ and 
$\max\left(\frac{e^{-\frac{5}{4}}}{4}, \frac{1-\sqrt{\frac{5}{8}}}{2} \right) = \frac{1-\sqrt{\frac{5}{8}}}{2}$.

\subsection{Proof of Corollary~\ref{corollary:covariance-estimation}} \label{sec:proof-corollary}
We use the same argument as in the beginning of the proof of Theorem~\ref{thm:main-stmt} in Section~\ref{sec:proof-main-stmt} but adjust the setting in which we apply Theorem~\ref{theorem:le-cam}. Specifically, we now let $\Theta  = \{\theta_{a} := C_{X_a} \,:\, X_a \sim \P_a,\, a \in \mathcal A \}$ with $C_X$ defined as in \eqref{eq:c-x} be the set of covariance operators, use the metric $(x,y) \mapsto \norm{x-y}{\H_k}$ for
$x,y\in\H_k$, and keep the remaining part of the setup the same. 
Hence, it remains to lower
bound $\norm{C_{X_{\theta_1}}-C_{X_{\theta_0}}}{\H_k}$. By using that HSIC is the RKHS norm of the cross-covariance operator, we obtain that
\begin{align}
  \norm{C_{X_{\theta_1}}-C_{X_{\theta_0}}}{\H_k} \stackrel{(i)}{\ge} 
  \Big|\underbrace{\norm{C_{X_{\theta_1}}}{\H_k}}_{=\HSIC_k\left(\P_{\theta_1}\right)}-\underbrace{\norm{C_{X_{\theta_0}}}{\H_k}}_{=\HSIC_k\left(\P_{\theta_0}\right)}\Big| = \left|F(\theta_1) - F(\theta_0)\right| \stackrel{(ii)}{\ge} 2s = \frac{2c}{\sqrt n}, \label{eq:corollary}
\end{align}
where $(i)$ holds by the reverse triangle inequality, $F$  is defined as in Section~\ref{sec:proof-main-stmt}, and $(ii)$ is guaranteed by~\eqref{eq:ti-bound} for $c>0$. We conclude as in the proof of Theorem~\ref{thm:main-stmt}(ii) to obtain the stated result.

\begin{ack}
This work was supported by the German Research Foundation (DFG) Research Training Group GRK 2153: Energy Status Data — Informatics Methods for its Collection, Analysis and Exploitation, and by the pilot program Core-Informatics of the Helmholtz Association (HGF).
\end{ack}

\bibliography{bib/collected_Zoltan.bib,bib/collected_plus.bib,bib/publications.bib}

\begin{thebibliography}{59}
\providecommand{\natexlab}[1]{#1}
\providecommand{\url}[1]{\texttt{#1}}
\expandafter\ifx\csname urlstyle\endcsname\relax
  \providecommand{\doi}[1]{doi: #1}\else
  \providecommand{\doi}{doi: \begingroup \urlstyle{rm}\Url}\fi

\bibitem[Albert et~al.(2022)Albert, Laurent, Marrel, and Meynaoui]{albert22adaptive}
M\'{e}lisande Albert, B\'{e}atrice Laurent, Amandine Marrel, and Anouar Meynaoui.
\newblock Adaptive test of independence based on {HSIC} measures.
\newblock \emph{The Annals of Statistics}, 50\penalty0 (2):\penalty0 858--879, 2022.

\bibitem[Anaya-Isaza and Mera-Jim{\'e}nez(2022)]{anaya22data}
Andr{\'e}s Anaya-Isaza and Leonel Mera-Jim{\'e}nez.
\newblock Data augmentation and transfer learning for brain tumor detection in magnetic resonance imaging.
\newblock \emph{IEEE Access}, 10:\penalty0 23217--23233, 2022.

\bibitem[Aronszajn(1950)]{aronszajn50theory}
Nachman Aronszajn.
\newblock Theory of reproducing kernels.
\newblock \emph{Transactions of the American Mathematical Society}, 68:\penalty0 337--404, 1950.

\bibitem[Berlinet and Thomas-Agnan(2004)]{berlinet04reproducing}
Alain Berlinet and Christine Thomas-Agnan.
\newblock \emph{Reproducing Kernel Hilbert Spaces in Probability and Statistics}.
\newblock Kluwer, 2004.

\bibitem[Bilodeau and Nangue(2017)]{bilodeau17tests}
Martin Bilodeau and Aur{\'e}lien~Guetsop Nangue.
\newblock Tests of mutual or serial independence of random vectors with applications.
\newblock \emph{Journal of Machine Learning Research}, 18:\penalty0 1--40, 2017.

\bibitem[Bouche et~al.(2023)Bouche, Flamary, d'Alch{\'e} Buc, Plougonven, Clausel, Badosa, and Drobinski]{bouche23wind}
Dimitri Bouche, R{\'e}mi Flamary, Florence d'Alch{\'e} Buc, Riwal Plougonven, Marianne Clausel, Jordi Badosa, and Philippe Drobinski.
\newblock Wind power predictions from nowcasts to 4-hour forecasts: a learning approach with variable selection.
\newblock \emph{Renewable Energy}, 211:\penalty0 938--947, 2023.

\bibitem[Camps-Valls et~al.(2010)Camps-Valls, Mooij, and Sch{\"{o}}lkopf]{camps10remote}
Gustavo Camps-Valls, Joris~M. Mooij, and Bernhard Sch{\"{o}}lkopf.
\newblock Remote sensing feature selection by kernel dependence measures.
\newblock \emph{{IEEE} {G}eoscience and {R}emote {S}ensing {L}etters}, 7\penalty0 (3):\penalty0 587--591, 2010.

\bibitem[Carmeli et~al.(2010)Carmeli, Vito, Toigo, and Umanit{\'a}]{carmeli10vector}
Claudio Carmeli, Ernesto~De Vito, Alessandro Toigo, and Veronica Umanit{\'a}.
\newblock Vector valued reproducing kernel {H}ilbert spaces and universality.
\newblock \emph{Analysis and Applications}, 8:\penalty0 19--61, 2010.

\bibitem[Chakraborty and Zhang(2019)]{chakraborty19distance}
Shubhadeep Chakraborty and Xianyang Zhang.
\newblock Distance metrics for measuring joint dependence with application to causal inference.
\newblock \emph{Journal of the American Statistical Association}, 114\penalty0 (528):\penalty0 1638--1650, 2019.

\bibitem[Climente-Gonz{\'a}lez et~al.(2019)Climente-Gonz{\'a}lez, Azencott, Kaski, and Yamada]{gonzalez19block}
Héctor Climente-Gonz{\'a}lez, Chlo{\'e}-Agathe Azencott, Samuel Kaski, and Makoto Yamada.
\newblock Block {HSIC} {L}asso: model-free biomarker detection for ultra-high dimensional data.
\newblock \emph{Bioinformatics}, 35\penalty0 (14):\penalty0 i427--i435, 2019.

\bibitem[Cohn(2013)]{cohn13measure}
Donald~L. Cohn.
\newblock \emph{Measure Theory}.
\newblock Birkh\"{a}user/Springer, second edition, 2013.

\bibitem[Diestel and Uhl(1977)]{diestel77vector}
Joseph Diestel and John~Jerry Uhl.
\newblock \emph{Vector Measures}.
\newblock American Mathematical Society. Providence, 1977.

\bibitem[Duchi(2007)]{duchi07derivations}
John Duchi.
\newblock Derivations for linear algebra and optimization.
\newblock \emph{Berkeley, California}, 3\penalty0 (1):\penalty0 2325--5870, 2007.

\bibitem[Fedorov et~al.(2024)Fedorov, Geenjaar, Wu, Sylvain, DeRamus, Luck, Misiura, Mittapalle, Hjelm, Plis, et~al.]{fedorov24multimodal}
Alex Fedorov, Eloy Geenjaar, Lei Wu, Tristan Sylvain, Thomas~P DeRamus, Margaux Luck, Maria Misiura, Girish Mittapalle, R~Devon Hjelm, Sergey~M Plis, et~al.
\newblock Self-supervised multimodal learning for group inferences from {MRI} data: Discovering disorder-relevant brain regions and multimodal links.
\newblock \emph{NeuroImage}, 285:\penalty0 120485, 2024.

\bibitem[Fellmann et~al.(2024)Fellmann, Blanchet-Scalliet, Helbert, Spagnol, and Sinoquet]{fellmann24kernel}
No{\'e} Fellmann, Christophette Blanchet-Scalliet, C{\'e}line Helbert, Adrien Spagnol, and Delphine Sinoquet.
\newblock Kernel-based sensitivity analysis for (excursion) sets.
\newblock \emph{Technometrics}, 2024.

\bibitem[Freitas~Gustavo et~al.(2023)Freitas~Gustavo, Hellström, and Verstraelen]{freitas23sensitivity}
Michael Freitas~Gustavo, Matti Hellström, and Toon Verstraelen.
\newblock Sensitivity analysis for {ReaxFF} reparametrization using the {H}ilbert--{S}chmidt independence criterion.
\newblock \emph{Journal of Chemical Theory and Computation}, 19\penalty0 (9):\penalty0 2557--2573, 2023.

\bibitem[Fukumizu et~al.(2008)Fukumizu, Gretton, Sun, and Sch{\"o}lkopf]{fukumizu08kernel}
Kenji Fukumizu, Arthur Gretton, Xiaohai Sun, and Bernhard Sch{\"o}lkopf.
\newblock Kernel measures of conditional dependence.
\newblock In \emph{Advances in Neural Information Processing Systems (NIPS)}, pages 498--496, 2008.

\bibitem[G\'orecki et~al.(2018)G\'orecki, Krzy\'sko, and Woly\'nski]{goreczki18independence}
Tomasz G\'orecki, Miroslaw Krzy\'sko, and Waldemar Woly\'nski.
\newblock Independence test and canonical correlation analysis based on the alignment between kernel matrices for multivariate functional data.
\newblock \emph{Artificial Intelligence Review}, pages 1--25, 2018.

\bibitem[Gretton et~al.(2005)Gretton, Bousquet, Smola, and Sch{\"o}lkopf]{gretton05measuring}
Arthur Gretton, Olivier Bousquet, Alex Smola, and Bernhard Sch{\"o}lkopf.
\newblock Measuring statistical dependence with {H}ilbert-{S}chmidt norms.
\newblock In \emph{Algorithmic Learning Theory (ALT)}, pages 63--78, 2005.

\bibitem[Gretton et~al.(2008)Gretton, Fukumizu, Teo, Song, Sch{\"o}lkopf, and Smola]{gretton08kernel}
Arthur Gretton, Kenji Fukumizu, Choon~Hui Teo, Le~Song, Bernhard Sch{\"o}lkopf, and Alexander Smola.
\newblock A kernel statistical test of independence.
\newblock In \emph{Advances in Neural Information Processing Systems (NIPS)}, pages 585--592, 2008.

\bibitem[Gretton et~al.(2012)Gretton, Borgwardt, Rasch, Sch{\"o}lkopf, and Smola]{gretton12kernel}
Arthur Gretton, Karsten Borgwardt, Malte Rasch, Bernhard Sch{\"o}lkopf, and Alexander Smola.
\newblock A kernel two-sample test.
\newblock \emph{Journal of Machine Learning Research}, 13\penalty0 (25):\penalty0 723--773, 2012.

\bibitem[Herrando-P{\'e}rez and Saltr{\'e}(2024)]{herrando24sensitivity}
Salvador Herrando-P{\'e}rez and Fr{\'e}d{\'e}rik Saltr{\'e}.
\newblock Estimating extinction time using radiocarbon dates.
\newblock \emph{Quaternary Geochronology}, 79:\penalty0 101489, 2024.

\bibitem[Kalinke and Szab{\'o}(2023)]{kalinke23nystrom}
Florian Kalinke and Zolt{\'a}n Szab{\'o}.
\newblock Nystr{\"o}m {M}-{H}ilbert-{S}chmidt independence criterion.
\newblock In \emph{Conference on Uncertainty in Artificial Intelligence (UAI)}, pages 1005--1015, 2023.

\bibitem[Le~Cam(1973)]{lecam73convergence}
Lucien Le~Cam.
\newblock Convergence of estimates under dimensionality restrictions.
\newblock \emph{The Annals of Statistics}, 1:\penalty0 38--53, 1973.

\bibitem[Lyons(2013)]{lyons13distance}
Russell Lyons.
\newblock Distance covariance in metric spaces.
\newblock \emph{The Annals of Probability}, 41:\penalty0 3284--3305, 2013.

\bibitem[Micchelli et~al.(2006)Micchelli, Xu, and Zhang]{micchelli06universal}
Charles Micchelli, Yuesheng Xu, and Haizhang Zhang.
\newblock Universal kernels.
\newblock \emph{Journal of Machine Learning Research}, 7:\penalty0 2651--2667, 2006.

\bibitem[Mooij et~al.(2016)Mooij, Peters, Janzing, Zscheischler, and Sch{\"o}lkopf]{mooij16distinguishing}
Joris Mooij, Jonas Peters, Dominik Janzing, Jakob Zscheischler, and Bernhard Sch{\"o}lkopf.
\newblock Distinguishing cause from effect using observational data: Methods and benchmarks.
\newblock \emph{Journal of Machine Learning Research}, 17:\penalty0 1--102, 2016.

\bibitem[Muandet et~al.(2011)Muandet, Fukumizu, Dinuzzo, and Sch{\"o}lkopf]{muandet12learning}
Krikamol Muandet, Kenji Fukumizu, Francesco Dinuzzo, and Bernhard Sch{\"o}lkopf.
\newblock Learning from distributions via support measure machines.
\newblock In \emph{Advances in Neural Information Processing Systems (NIPS)}, pages 10--18, 2011.

\bibitem[M{\"u}ller(1997)]{muller97integral}
Alfred M{\"u}ller.
\newblock Integral probability metrics and their generating classes of functions.
\newblock \emph{Advances in Applied Probability}, 29:\penalty0 429--443, 1997.

\bibitem[Pfister et~al.(2018)Pfister, B{\"u}hlmann, Sch{\"o}lkopf, and Peters]{pfister18kernel}
Niklas Pfister, Peter B{\"u}hlmann, Bernhard Sch{\"o}lkopf, and Jonas Peters.
\newblock Kernel-based tests for joint independence.
\newblock \emph{Journal of the Royal Statistical Society: Series B (Statistical Methodology)}, 80\penalty0 (1):\penalty0 5--31, 2018.

\bibitem[Podkopaev et~al.(2023)Podkopaev, Bl\"{o}baum, Kasiviswanathan, and Ramdas]{podkopaev23sequential}
Aleksandr Podkopaev, Patrick Bl\"{o}baum, Shiva Kasiviswanathan, and Aaditya Ramdas.
\newblock Sequential kernelized independence testing.
\newblock In \emph{International Conference on Machine Learning ({ICML})}, pages 27957--27993, 2023.

\bibitem[Quadrianto et~al.(2009)Quadrianto, Song, and Smola]{quadrianto09kernelized}
Novi Quadrianto, Le~Song, and Alex Smola.
\newblock Kernelized sorting.
\newblock In \emph{Advances in Neural Information Processing Systems (NIPS)}, pages 1289--1296, 2009.

\bibitem[Saitoh and Sawano(2016)]{saitoh16theory}
Saburou Saitoh and Yoshihiro Sawano.
\newblock \emph{Theory of Reproducing Kernels and Applications}.
\newblock Springer Singapore, 2016.

\bibitem[Sch{\"{o}}lkopf et~al.(2021)Sch{\"{o}}lkopf, Locatello, Bauer, Ke, Kalchbrenner, Goyal, and Bengio]{scholkopf21causal}
Bernhard Sch{\"{o}}lkopf, Francesco Locatello, Stefan Bauer, Nan~Rosemary Ke, Nal Kalchbrenner, Anirudh Goyal, and Yoshua Bengio.
\newblock Toward causal representation learning.
\newblock \emph{Proceedings of the {IEEE}}, 109\penalty0 (5):\penalty0 612--634, 2021.

\bibitem[Sejdinovic et~al.(2013{\natexlab{a}})Sejdinovic, Gretton, and Bergsma]{sejdinovic13kernel}
Dino Sejdinovic, Arthur Gretton, and Wicher Bergsma.
\newblock A kernel test for three-variable interactions.
\newblock In \emph{Advances in Neural Information Processing Systems (NIPS)}, pages 1124--1132, 2013{\natexlab{a}}.

\bibitem[Sejdinovic et~al.(2013{\natexlab{b}})Sejdinovic, Sriperumbudur, Gretton, and Fukumizu]{sejdinovic13equivalence}
Dino Sejdinovic, Bharath Sriperumbudur, Arthur Gretton, and Kenji Fukumizu.
\newblock Equivalence of distance-based and {RKHS}-based statistics in hypothesis testing.
\newblock \emph{Annals of Statistics}, 41:\penalty0 2263--2291, 2013{\natexlab{b}}.

\bibitem[Shekhar et~al.(2023)Shekhar, Kim, and Ramdas]{shekhar23permutation}
Shubhanshu Shekhar, Ilmun Kim, and Aaditya Ramdas.
\newblock A permutation-free kernel independence test.
\newblock \emph{Journal of Machine Learning Research}, 24\penalty0 (369):\penalty0 1--68, 2023.

\bibitem[Sheng and Sriperumbudur(2023)]{sheng23distance}
Tianhong Sheng and Bharath~K. Sriperumbudur.
\newblock On distance and kernel measures of conditional independence.
\newblock \emph{Journal of Machine Learning Research}, 24\penalty0 (7):\penalty0 1--16, 2023.

\bibitem[Smola et~al.(2007)Smola, Gretton, Song, and Sch{\"o}lkopf]{smola07hilbert}
Alexander Smola, Arthur Gretton, Le~Song, and Bernhard Sch{\"o}lkopf.
\newblock A {H}ilbert space embedding for distributions.
\newblock In \emph{Algorithmic Learning Theory (ALT)}, pages 13--31, 2007.

\bibitem[Song et~al.(2007)Song, Smola, Gretton, and Borgwardt]{song07dependence}
Le~Song, Alexander~J. Smola, Arthur Gretton, and Karsten~M. Borgwardt.
\newblock A dependence maximization view of clustering.
\newblock In \emph{International Conference on Machine Learning (ICML)}, pages 815--822, 2007.

\bibitem[Song et~al.(2012)Song, Smola, Gretton, Bedo, and Borgwardt]{song12feature}
Le~Song, Alex Smola, Arthur Gretton, Justin Bedo, and Karsten Borgwardt.
\newblock Feature selection via dependence maximization.
\newblock \emph{Journal of Machine Learning Research}, 13\penalty0 (1):\penalty0 1393--1434, 2012.

\bibitem[Sriperumbudur et~al.(2010)Sriperumbudur, Gretton, Fukumizu, Sch{\"o}lkopf, and Lanckriet]{sriperumbudur10hilbert}
Bharath Sriperumbudur, Arthur Gretton, Kenji Fukumizu, Bernhard Sch{\"o}lkopf, and Gert Lanckriet.
\newblock Hilbert space embeddings and metrics on probability measures.
\newblock \emph{Journal of Machine Learning Research}, 11:\penalty0 1517--1561, 2010.

\bibitem[Sriperumbudur et~al.(2011)Sriperumbudur, Fukumizu, and Lanckriet]{sriperumbudur11universality}
Bharath Sriperumbudur, Kenji Fukumizu, and Gert Lanckriet.
\newblock Universality, characteristic kernels and {RKHS} embedding of measures.
\newblock \emph{Journal of Machine Learning Research}, 12:\penalty0 2389--2410, 2011.

\bibitem[Steinwart(2001)]{steinwart01influence}
Ingo Steinwart.
\newblock On the influence of the kernel on the consistency of support vector machines.
\newblock \emph{Journal of Machine Learning Research}, 6\penalty0 (3):\penalty0 67--93, 2001.

\bibitem[Steinwart and Christmann(2008)]{steinwart08support}
Ingo Steinwart and Andreas Christmann.
\newblock \emph{Support Vector Machines}.
\newblock Springer, 2008.

\bibitem[Stenger et~al.(2020)Stenger, Gamboa, Keller, and Iooss]{stenger2020optimal}
Jerome Stenger, Fabrice Gamboa, Merlin Keller, and Bertrand Iooss.
\newblock Optimal uncertainty quantification of a risk measurement from a thermal-hydraulic code using canonical moments.
\newblock \emph{International Journal for Uncertainty Quantification}, 10\penalty0 (1), 2020.

\bibitem[Szab{\'o} and Sriperumbudur(2018)]{szabo18characteristic2}
Zolt{\'a}n Szab{\'o} and Bharath~K. Sriperumbudur.
\newblock Characteristic and universal tensor product kernels.
\newblock \emph{Journal of Machine Learning Research}, 18\penalty0 (233):\penalty0 1--29, 2018.

\bibitem[Sz{\'e}kely and Rizzo(2009)]{szekely09brownian}
G{\'a}bor~J. Sz{\'e}kely and Maria~L. Rizzo.
\newblock Brownian distance covariance.
\newblock \emph{The Annals of Applied Statistics}, 3:\penalty0 1236--1265, 2009.

\bibitem[Sz{\'e}kely et~al.(2007)Sz{\'e}kely, Rizzo, and Bakirov]{szekely07measuring}
G{\'a}bor~J. Sz{\'e}kely, Maria~L. Rizzo, and Nail~K. Bakirov.
\newblock Measuring and testing dependence by correlation of distances.
\newblock \emph{The Annals of Statistics}, 35:\penalty0 2769--2794, 2007.

\bibitem[Tolstikhin et~al.(2016)Tolstikhin, Sriperumbudur, and Sch{\"o}lkopf]{tolstikhin16minimax}
Ilya Tolstikhin, Bharath Sriperumbudur, and Bernhard Sch{\"o}lkopf.
\newblock Minimax estimation of maximal mean discrepancy with radial kernels.
\newblock In \emph{Advances in Neural Information Processing Systems (NIPS)}, pages 1930--1938, 2016.

\bibitem[Tolstikhin et~al.(2017)Tolstikhin, Sriperumbudur, and Muandet]{tolstikhin17minimax}
Ilya Tolstikhin, Bharath Sriperumbudur, and Krikamol Muandet.
\newblock Minimax estimation of kernel mean embeddings.
\newblock \emph{Journal of Machine Learning Research}, 18:\penalty0 1--47, 2017.

\bibitem[Tsybakov(2009)]{tsybakov09nonparametric}
Alexandre~B. Tsybakov.
\newblock \emph{Introduction to Nonparametric Estimation}.
\newblock Springer, 2009.

\bibitem[Veiga(2015)]{veiga15global}
Sebastien~De Veiga.
\newblock Global sensitivity analysis with dependence measures.
\newblock \emph{Journal of Statistical Computation and Simulation}, 85\penalty0 (7):\penalty0 1283--1305, 2015.

\bibitem[Wang et~al.(2022)Wang, Du, Zhang, and Shi]{wang22rank}
Andi Wang, Juan Du, Xi~Zhang, and Jianjun Shi.
\newblock Ranking features to promote diversity: An approach based on sparse distance correlation.
\newblock \emph{Technometrics}, 64\penalty0 (3):\penalty0 384--395, 2022.

\bibitem[Wehbe and Ramdas(2015)]{wehbe15nonparametric}
Leila Wehbe and Aaditya Ramdas.
\newblock Nonparametric independence testing for small sample sizes.
\newblock In \emph{International Joint Conference on Artificial Intelligence ({IJCAI})}, pages 3777--3783, 2015.

\bibitem[Wendland(2005)]{wendland05scattered}
Holger Wendland.
\newblock \emph{Scattered data approximation}.
\newblock Cambridge University Press, 2005.

\bibitem[Yamada et~al.(2014)Yamada, Jitkrittum, Sigal, Xing, and Sugiyama]{yamada14high}
Makoto Yamada, Wittawat Jitkrittum, Leonid Sigal, Eric~P. Xing, and Masashi Sugiyama.
\newblock High-dimensional feature selection by feature-wise kernelized {L}asso.
\newblock \emph{Neural Computation}, 26\penalty0 (1):\penalty0 185--207, 2014.

\bibitem[Zhou et~al.(2019)Zhou, Chen, and Huang]{zhou19covestimators}
Yang Zhou, Di-Rong Chen, and Wei Huang.
\newblock A class of optimal estimators for the covariance operator in reproducing kernel {H}ilbert spaces.
\newblock \emph{Journal of Multivariate Analysis}, 169:\penalty0 166--178, 2019.

\bibitem[Zolotarev(1983)]{zolotarev83probability}
V.~Zolotarev.
\newblock Probability metrics.
\newblock \emph{Theory of Probability and its Applications}, 28:\penalty0 278--302, 1983.

\end{thebibliography}

\newpage
\appendix
\setcounter{equation}{0}
\renewcommand{\theequation}{\thesection.\arabic{equation}}

\section{Auxiliary Result}\label{sec:auxiliary-results}

In this section, we collect an auxiliary result. Lemma~\ref{lemma:bound-kl}
presents an upper bound on the Kullback-Leibler divergence between multivariate
normal distributions.

\begin{lemmaA}[Upper bound on KL divergence]\label{lemma:bound-kl} Let $d=\sum_{m=1}^Md_m$, with $d_m\in \mathbb{N}_{>0}$ ($m\in [M]$). Fix $i \in
  [d_1]$. Let $j=i+1$, $\P_{\theta_0} = \mathcal N (\bm 0_d,\b I_d)$, and
  $\P_{\theta_1} = \mathcal N(\bm \mu_1,\bm \Sigma_1)$, with $\bm\mu_1 =
  \frac{1}{\sqrt dn}\bm1_d \in \R^d$, and $\bm \Sigma_1 = \bm \Sigma(i,j,\rho_n)
  \in \R^{d\times d}$ defined as in \eqref{eq:Sigma-rho}
  ($\rho_n
  \in (0,1)$). Then, for $2\le n \in \mathbb N$,
  \begin{align*}
    \mathrm{KL}(\P_{\theta_1}^n||\P_{\theta_0}^n) \le \frac{1}{2n}+ \frac n 2 \frac{\rho^2_n}{1-\rho^2_n}.
  \end{align*}
  In particular, for $\rho_n^2 = 1/n$, it holds that $\mathrm{KL}(\P_{\theta_1}^n||\P_{\theta_0}^n) \le \frac{5}{4}$.

\begin{proof}
With $\bm{\mu}_0=\b 0_d$ and $\bm \Sigma_0 = \b I_d$, we obtain that
\begin{align*}
    \mathrm{KL}(\P_{\theta_1}^n||\P_{\theta_0}^n) &\stackrel{(a)}{=} \sum_{i\in[n]}\mathrm{KL}(\P_{\theta_1}||\P_{\theta_0}) \\
    & \stackrel{(b)}{=} \frac{n}{2}\left[
    \trace(\b\Sigma_0^{-1} \b\Sigma_1) + (\bm\mu_0-\bm\mu_1)\T\b\Sigma_0^{-1}(\bm\mu_0-\bm\mu_1) - d +\ln\left(\frac{\left|\b\Sigma_0\right|}{\left|\b\Sigma_1\right|}\right)
    \right] \\
    &= \frac{n}{2}\Bigg[ \underbrace{\trace(\b\Sigma_1)}_{=d} + \underbrace{\left\|\bm \mu_1\right\|_{\R^d}^2}_{=\frac{1}{n^2}} - d + \ln\Bigg(\underbrace{\frac{1}{\left|\bm \Sigma_1\right|}}_{\stackrel{(c)}{=}\frac{1}{1-\rho_n^2}}\Bigg)\Bigg] \\
    &= \frac{1}{2n}+\frac{n}{2} \ln \left(\frac{1}{1-\rho_n^2}\right) \stackrel{(d)}{\le}\frac{1}{2n} + \frac n 2 \frac{\rho^2_n}{1-\rho^2_n} \stackrel{(e)}{\le}   \frac54,
\end{align*}
where (a) is implied by  Lemma~\ref{lemma:kl-sum}, (b) follows from Lemma~\ref{lemma:kl-gaussians},  (c) follows from the definition of the determinant, (d)
is the consequence of the inequality $\ln(x) \le x-1$ holding for $x> 0$, and (e) holds for $n\ge2$ and $\rho_n^2 = 1/n$ as
\begin{align*}
   \frac n 2 \underbrace{\frac{1/n}{1-1/n}}_{  \frac{1}{n-1} } \le 1  \iff \frac n 2   \frac{1}{n-1} \le  1 \iff n \le 2(n-1) \iff n\ge 2 ,
\end{align*}
and in this case (for $n\ge2$) one has that $\frac {1}{2n}\le \frac14$.
\end{proof}
\end{lemmaA}

\section{External Theorems}
\label{sec:external-theorems}

For self-completeness, we include the external statements that we use. The well-known result by Bochner, stated in Theorem~\ref{thm:bochner}, completely characterizes continuous bounded translation-invariant kernels. 
Theorem~\ref{thm:bharath-char-func}
allows expressing MMD with  continuous bounded translation-invariant kernels in terms of
characteristic functions, and Theorem~\ref{thm:bharath-full-Rd} gives an equivalent condition for a continuous bounded translation-invariant kernel to be characteristic. Theorem~\ref{thm:szabo-thm4} connects characteristic kernels to characteristic product kernels and to $\mathcal{I}$-characteristic product kernels on $\R^d$ (we include only the part relevant to our paper for brevity). 
We recall Le Cam's method in
Theorem~\ref{theorem:le-cam} and collect results on the Kullback-Leibler
divergence in Lemma~\ref{lemma:kl-sum} and Lemma~\ref{lemma:kl-gaussians}.

\begin{theoremA}[Bochner; Theorem~6.6; \citet{wendland05scattered}] \label{thm:bochner}
A continuous function $\kappa : \R^d \to \R$ is positive definite if and only if it is the Fourier transform of a finite nonnegative Borel measure $\Lambda$ on $\R^d$, that is,
\begin{align*} 
 \kappa(\b x) = \int_{\R^d}e^{-i\ip{\b x,\bm \omega}{}}\d\Lambda(\bm \omega),\quad \text{for all } \b x \in \R^d.
\end{align*}
\end{theoremA}

\begin{theoremA}[Corollary 4(i); \citet{sriperumbudur10hilbert}] \label{thm:bharath-char-func} Let $k:\R^d \times \R^d \to \R$ be a  continuous bounded translation-invariant kernel. Then, for any $\P,\Q \in \mathcal M_1^+\left(\R^d\right)$,
  \begin{align*}
  \MMD_k^2(\P,\Q) = \norm{\psi_\P-\psi_\Q}{L^2\left(\R^d,\Lambda_k\right)}^2,
  \end{align*}
  with $\psi_\P$ and $\psi_\Q$ being the characteristic functions of $\P$ and $\Q$, respectively, and $\Lambda_k$ defined in \eqref{eq:bochner}.
  
\end{theoremA}

\begin{theoremA}[Theorem 9; \citet{sriperumbudur10hilbert}] \label{thm:bharath-full-Rd} Suppose $k:\R^d \times \R^d \to \R$ is a  continuous bounded translation-invariant kernel. Then $k$ is characteristic if and only if $\operatorname{supp}(\Lambda_k) = \R^d$, with $\Lambda_k$ defined as in \eqref{eq:bochner}.
\end{theoremA}

\begin{theoremA}[Theorem 4; \citet{szabo18characteristic2}]
\label{thm:szabo-thm4} Suppose $k_m : \R^{d_m} \times \R^{d_m} \to \R$ is
continuous bounded and translation-invariant kernel for all $m\in[M]$. Then the
following statements are equivalent:
  \begin{enumerate}[itemindent=1.7em,leftmargin=0pt,labelsep=4pt,itemsep=4pt,topsep=-3pt]
  \item[(i)] $(k_m)_{m=1}^M$-s are characteristic;
    \item[(ii)] $\otimes_{m=1}^Mk_m$ is characteristic;
    \item[(iii)] $\otimes_{m=1}^Mk_m$ is $\mathcal{I}$-characteristic.
  \end{enumerate}
\end{theoremA}

The next statement follows directly from \citet[Eq.~(2.9)]{tsybakov09nonparametric} and \citet[Theorem~2.2]{tsybakov09nonparametric}.

\begin{theoremA}[Theorem 2.2; \citet{tsybakov09nonparametric}]
\label{theorem:le-cam}
Let $\X$ be a measurable space, $(\Theta,d)$ is a semi-metric space, and $\mathcal P_{\Theta} =
\{\P_\theta : \theta \in  \Theta\}$ is a class of probability measures
on $\X$ indexed by $\Theta$. We
observe data $D \sim \P_{\theta} \in \mathcal P_{ \Theta}$ with some unknown
parameter $\theta$. The goal is to estimate $\theta$. Let $\hat \theta = \hat \theta(D)$ be an estimator of $\theta$ based on~$D$.
Assume that there exist $\theta_0,\theta_1\in \Theta$ such that
$d(\theta_0,\theta_1) \ge 2s > 0$ and $\mathrm{KL}(\P_{\theta_1}||\P_{\theta_0})
\le \alpha <\infty$ for $\alpha > 0$. Then
\begin{align*}
  \inf_{\hat \theta}\sup_{\theta\in \Theta}\P_\theta\left(d\left(\hat \theta,\theta\right) \ge s\right) \ge \max \left(\frac{e^{-\alpha}}{4},\frac{1-\sqrt{\alpha/2}}{2}\right).
\end{align*}

\end{theoremA}

We have the following property of the Kullback-Leibler divergence for product measures \citep[p.~85]{tsybakov09nonparametric}.
\begin{lemmaA}[KL divergence of product measures]
\label{lemma:kl-sum}
Let $\P=\otimes_{i=1}^n\P_i$ and $\Q=\otimes_{i=1}^n\Q_i$. Then
\begin{align*}
    \mathrm{KL}(\P||\Q) = \sum_{i\in[n]}\mathrm{KL}(\P_i||\Q_i).
\end{align*}
\end{lemmaA}

The following lemma \cite[p.~13]{duchi07derivations} shows that the Kullback-Leibler divergence of multivariate Gaussians can be computed in closed form.

\begin{lemmaA}[KL divergence of Gaussians]
\label{lemma:kl-gaussians}
The KL divergence of two normal distributions $\mathcal{N}(\bm
\mu_1,\b\Sigma_1)$ and $\mathcal{N}(\bm \mu_0,\b\Sigma_0)$ on $\R^d$   is
\begin{align*}
  \mathrm{KL}(\mathcal{N}(\bm \mu_1,\b\Sigma_1)||\mathcal{N}(\bm \mu_0,\b\Sigma_0)) = \frac{
    \trace(\b\Sigma_0^{-1} \b\Sigma_1) + (\bm\mu_0-\bm\mu_1)\T\b\Sigma_0^{-1}(\bm\mu_0-\bm\mu_1) - d +\ln\left(\frac{\left|\b\Sigma_0\right|}{\left|\b\Sigma_1\right|}\right)
    }{2}.
\end{align*}
\end{lemmaA}

\end{document}